\documentstyle[12pt,amssymb,amscd]{amsart}  

\theoremstyle{remark}

\theoremstyle{definition}

\numberwithin{equation}{subsection}

\newcommand{\bbC}{{\Bbb C}}

\newcommand{\cF}{{\cal F}}

\newcommand{\cO}{{\cal O}}

\newcommand{\cL}{{\cal L}}

\newcommand{\cV}{{\cal V}}
\newcommand{\cA}{{\cal A}}

\newcommand{\cC}{{\cal C}}
\newcommand{\cE}{{\cal E}}

\newcommand{\cU}{{\cal U}}

\newcommand{\cS}{{\cal S}}
\newcommand{\cT}{{\cal T}}
\newcommand{\cQ}{{\cal Q}}

\newcommand{\isomo}{\overset{\sim}{=}}

\newcommand{\id}{\operatorname{id}}
\newcommand{\Tr}{\operatorname{Tr}}
\newcommand{\shHom}{\underline{\operatorname{Hom}}}
\newcommand{\shEnd}{\underline{\operatorname{End}}}
\newcommand{\shExt}{\underline{\operatorname{Ext}}}

\newcommand{\Ext}{\operatorname{Ext}}

\newcommand{\ch}{\operatorname{ch}}

\newcommand{\pr}{\operatorname{pr}}

\newcommand{\ECA}{{\mathcal ECA}}

\newcommand{\EVA}{{\mathcal EVA}}
\newcommand{\CEXT}{{\mathcal CE\!\scriptstyle{XT}}}
\newcommand{\VEXT}{{\mathcal VE\!\scriptstyle{XT}}}
\newcommand{\DR}{{\Omega^\bullet_X}}

\newcommand\dual[1]{{#1}^{\vee}}
\newcommand{\ip}{{\langle\ ,\ \rangle}}
\newcommand{\Atiyah}{\alpha}
\newcommand{\pont}{{\scriptstyle{\Pi}}}

\begin{document}

\title{Vertex Algebroids II}

\author{Paul Bressler}
\address{}
\email{bressler@@math.arizona.edu}
\date{\today}
\maketitle

\section{Introduction}
In this note we determine the obstruction to triviality  of the stack of
exact vertex algebroids thereby recovering the result of \cite{GMS}.

The stack $\EVA_{\cO_X}$ of exact vertex $\cO_X$-algebroids is a torsor under
the stack in Picard groupoids $\ECA_{\cO_X}$ of exact Courant
$\cO_X$-algebroids. The latter is equivalent to the stack of torsors under
$\Omega^2_X @>>> \Omega^{3,cl}$. Therefore, $\ECA_{\cO_X}$-torsors are classified
by $H^2(X;\Omega^2_X @>>> \Omega^{3,cl})$. The goal of the present note is to
determine the class of $\EVA_{\cO_X}$.

The first step toward this goal is to replace $\EVA_{\cO_X}$ by the equivalent
$\ECA_{\cO_X}$-torsor $\CEXT_{\cO_X}(\cA_{\Omega^1_X})_\ip$,
whose (locally defined) objects are certain Courant algebroids which are extensions by
$\Omega^1_X$ of the Lie $\cO_X$-algebroid $\cA_{\Omega^1_X}$, the Atiyah algebra
of the sheaf $\Omega^1_X$. Any such extension induces an $\cA_{\Omega^1_X}$-invariant
symmetric pairing $\ip$ on the Lie algebra $\shEnd_{\cO_X}(\Omega^1_X)$. The objects of
$\CEXT_{\cO_X}(\cA_{\Omega^1_X})_\ip$ are those for which $\ip$ is given by
the {\em negative} of the trace of the product of endomorphisms.

We show that $\EVA_{\cO_X}$ and $\CEXT_{\cO_X}(\cA_{\Omega^1_X})_\ip$ are
anti-equivalent as $\ECA_{\cO_X}$-torsors by adapting
the strategy of \cite{BD} to the present setting and making use of
(the degree zero part of) the unique vertex $\DR$-algebroid
constructed in \cite{B}. It follows that the classes of $\EVA_{\cO_X}$ and
$\CEXT_{\cO_X}(\cA_{\Omega^1_X})_\ip$ in $H^2(X;\Omega^2_X @>>> \Omega^{3,cl})$
are negatives of each other. The advantage of passing to
$\CEXT_{\cO_X}(\cA_{\Omega^1_X})_{\Tr}$ has to do with the fact
that Courant algebroids are objects of ``classical'' nature. In
particular they are $\cO_X$-modules (as opposed to vertex
algebroids).

The second step is the determination of the class of
$\CEXT_{\cO_X}(\cA_{\Omega^1_X})_\ip$ which is achieved a the
more general framework. Namely, we consider a transitive Lie
$\cO_X$ algebroid $\cA$, locally free of finite rank over $\cO_X$,
and denote by $\frak g$ the kernel of the anchor map. Thus, $\frak
g$ is a sheaf of Lie algebras in $\cO_X$-modules. If $\widehat\cA$
is a Courant algebroid, such that the associated Lie algebroid is
identified with $\cA$, then $\widehat\cA$ is an extension of $\cA$
by $\Omega^1_X$. The symmetric bilinear form on $\widehat\cA$ induces
an $\cA$-invariant symmetric bilinear form on on $\frak g$.

Suppose $\cA$, $\frak g$ as above, and $\ip$ an $\cA$-invariant symmetric
bilinear form on on $\frak g$. Let $\CEXT_{\cO_X}(\cA)_\ip$ denote the stack of
whose (locally defined) objects are pairs, consisting of a Courant algebroid
$\widehat\cA$ together with an identification of the associated Lie algebroid
with $\cA$, such that
the symmetric bilinear form induced on $\frak g$ coincides with $\ip$. We show that,
if $\cA$ admits a flat connection locally on $X$, then
$\CEXT_{\cO_X}(\cA)_\ip$ has a natural structure of an $\ECA_{\cO_X}$-torsor
and calculate its characteristic class in $H^2(X;\Omega^2_X @>>> \Omega^{3,cl})$.
It turns out that it is equal to the Pontryagin class naturally associated to
the pair $(\cA,\ip)$. For example, the Pontryagin class of $(\cA_{\Omega^1_X},\ip)$
(where the symmetric pairing is given by the negative of the trace of
the product of endomorphisms) is equal to $-2\ch(\Omega^1_X)$. Hence, the class
of $\EVA_{\cO_X}$ in $H^2(X;\Omega^2_X @>>> \Omega^{3,cl})$ is equal to
$2\ch(\Omega^1_X)$.

The paper is organized as follows. In Section 2 we recall the basic material on
Courant algebroids. In Section 3 we study Courant extensions of
transitive Lie algebroids and classification thereof. In Section 4 we recall the
basic material on vertex algebroids and their relationship to Courant algebroids.
We explain how the construction of \cite{B} gives an example of a vertex extension
(of the Atiyah algebra of the cotangent bundle) and use this example to classify
exact vertex algebroids.

The fact that the Pontryagin class of a principal bundle (defined with respect to
an invariant symmetric pairing on the Lie algebra of the structure group) is the
obstruction to the existence of a Courant extension of the Atiyah algebra of the
principal bundle was pointed out by P.~\v Severa in \cite{S} together with the
constructions of \ref{subsection:CEXT-flat}.

\section{Courant algebroids}

\subsection{Courant algebroids}
A {\em Courant $\cO_X$-algebroid} is an $\cO_X$-module $\cQ$
equipped with
\begin{enumerate}
\item a structure of a Leibniz $\bbC$-algebra
\[
[\ ,\ ] : \cQ\otimes_\bbC\cQ @>>> \cQ \ ,
\]

\item
an $\cO_X$-linear map of Leibniz algebras (the anchor map)
\[
\pi : \cQ @>>> \cT_X \ ,
\]

\item
a symmetric $\cO_X$-bilinear pairing
\[
\ip : \cQ\otimes_{\cO_X}\cQ @>>> \cO_X \ ,
\]

\item
a derivation
\[
\partial : \cO_X @>>> \cQ
\]
such that $\pi\circ\partial = 0$
\end{enumerate}
which satisfy
\begin{eqnarray}
[q_1,fq_2] & = & f[q_1,q_2] + \pi(q_1)(f)q_2 \label{leibniz}\\
\langle [q,q_1],q_2\rangle + \langle q_1,[q,q_2]\rangle & = & \pi(q)(\langle q_1, q_2\rangle)
\label{ip-invar}\\
\left[q,\partial(f)\right] & = & \partial(\pi(q)(f)) \label{bracket-o-courant}\\
\langle q,\partial(f)\rangle & = & \pi(q)(f) \label{ip-o}\\
\left[q_1,q_2\right] + [q_2,q_1] & = & \partial(\langle q_1, q_2\rangle) \label{ip-symm}
\end{eqnarray}
for $f\in\cO_X$ and $q,q_1,q_2\in\cQ$.

\subsubsection{}
A morphism of Courant $\cO_X$-algebroids is an $\cO_X$-linear map
of Leibnitz algebras which commutes with the respective anchor
maps and derivations and preserves the respective pairings.

\subsection{The associated Lie algebroid}
Suppose that $\cQ$ is a Courant $\cO_X$-algebroid. Let
\begin{eqnarray*}
\Omega_\cQ & \stackrel{def}{=} & \cO_X\partial(\cO_X)\subset\cQ \ , \\
\overline\cQ & \stackrel{def}{=} & \cQ/\Omega_\cQ \ .
\end{eqnarray*}
Note that the symmetrization of the Leibniz bracket on $\cQ$ takes
values in $\Omega_\cQ$.

For $q\in\cQ$, $f,g\in\cO_X$
\begin{eqnarray*}
[q,f\partial(g)] & = & f[q,\partial(g)] + \pi(q)(f)\partial(g) \\
& = & f\partial(\pi(q)(g)) + \pi(q)(f)\partial(g)
\end{eqnarray*}
which shows that $[\cQ,\Omega_\cQ]\subseteq\Omega_\cQ$. Therefore,
the Leibniz bracket on $\cQ$ descends to the Lie bracket
\begin{equation}\label{LAbracket}
[\ ,\ ] : \overline\cQ\otimes_\bbC\overline\cQ @>>> \overline\cQ\ .
\end{equation}

Since $\pi$ is $\cO_X$-linear and $\pi\circ\partial = 0$, $\pi$
vanishes on $\Omega_\cQ$ and factors through the map
\begin{equation}\label{LAanchor}
\pi : \overline\cQ @>>> \cT_X \ .
\end{equation}

\subsubsection{Lemma}
The bracket \eqref{LAbracket} and the anchor \eqref{LAanchor}
determine the structure of a Lie $\cO_X$-algebroid on
$\overline\cQ$.

\subsection{Transitive Courant algebroids}
A Courant $\cO_X$-algebroid is called {\em transitive} if the anchor
is surjective.

\subsubsection{Remark}
A Courant $\cO_X$-algebroid $\cQ$ is transitive if and only if the
associated Lie $\cO_X$-algebroid is.

\subsubsection{}
Suppose  that  $\cQ$ is a transitive Courant $\cO_X$-algebroid. The
derivation $\partial$ induces the $\cO_X$-linear map
\[
i : \Omega^1_X @>>> \cQ \ .
\]
Since $\langle q, \alpha\rangle = \iota_{\pi(q)}\alpha$, it
follows that the map $i$ is adjoint to the anchor map $\pi$. The
surjectivity of the latter implies that $i$ is injective. Since,
in addition, $\pi\circ i = 0$ the sequence
\begin{equation}\label{ses:assoc-Lie-alg}
0 @>>> \Omega^1_X @>{i}>> \cQ @>>> \overline\cQ @>>> 0
\end{equation}
is exact. Moreover, $i$ is isotropic with respect to the
symmetric pairing.

\subsubsection{Definition}
A {\em connection} on a transitive Courant $\cO_X$-algebroid $\cQ$
is a $\cO_X$-linear isotropic section of the anchor map $\cQ @>>>
\cT_X$.

\subsubsection{Definition}
A {\em flat connection} on a transitive Courant $\cO_X$-algebroid
$\cQ$ is a $\cO_X$-linear section of the anchor map which is
morphism of Leibniz algebras.

\subsubsection{Remark}
As a consequence of \eqref{ip-symm} a flat connection is a connection.

\subsection{Exact Courant algebroids}
\subsubsection{Definition}
The Courant algebroid $\cQ$ is called {\em exact} if the anchor
map $\pi: \overline{\cQ} @>>> \cT_X$ is an isomorphism.

We denote the stack of exact Courant $\cO_X$-algebroids by
$\ECA_{\cO_X}$.

\subsubsection{}
A morphism of exact Courant algebroids induces a morphism of
respective extensions of $\cT_X$ by $\Omega^1_X$, hence is an
isomorphism of $\cO_X$-modules. It is clear that the inverse
isomorphism is a morphism of Courant $\cO_X$-algebroids.
Therefore, $\ECA_{\cO_X}$ is a stack in groupoids.

\subsubsection{}
The evident morphism
$\ECA_{\cO_X} @>>> \shExt^1_{\cO_X}(\cT_X,\Omega^1_X)$
is faithful. The natural structure of a $\bbC$-vector
space in categories on $\shExt^1_{\cO_X}(\cT_X,\Omega^1_X)$
restricts to one on $\ECA_{\cO_X}$. In particular,
$\ECA_{\cO_X}$ is a stack in Picard groupoids.

\subsubsection{Connections}
Suppose that $\cQ$ is an exact Courant $\cO_X$-algebroid. Let
$\cC(\cQ)$ denote the sheaf of (locally defined) connections on
$\cQ$.

\subsubsection{Lemma}
$\cC(\cQ)$ is an $\Omega^2_X$-torsor.
\begin{pf}
The difference of two sections of the anchor map $\cQ @>>> \cT_X$
is a map $\cT_X @>>> \Omega^1_X$ or, equivalently, a section of
$\Omega^1_X\otimes_{\cO_X}\Omega^1_X$. The difference of two
isotropic sections gives rise to a skew-symmetric tensor, i.e.
a section of $\Omega^2_X$.
\end{pf}

\subsubsection{Curvature}
For a (locally defined)
connection $\nabla$ the formula
\[
(\xi,\xi_1,\xi_2)\mapsto\iota_\xi([\nabla(\xi_1),\nabla(\xi_2)] -
\nabla([\xi_1,\xi_2])
\]
defines a section, denoted $c(\nabla)$, of
$\Omega^1_X\otimes_{\cO_X}\Omega^1_X\otimes_{\cO_X}\Omega^1_X$
called {\em the curvature} of the connection $\nabla$.

A connection $\nabla$ is called {\em flat} if $c(\nabla)=0$.

\subsubsection{Lemma}
\begin{enumerate}
\item The tensor $c(\nabla)$ is skew-symmetric, i.e. a
 section of $\Omega^3_X$.

\item The differential form $c(\nabla)$ is closed.

\item For $\alpha\in\Omega^2_X$, $c(\nabla+\alpha) =
c(\nabla) + d\alpha$.
\end{enumerate}

\subsubsection{Exact Courant algebroids with connection}
Pairs $(\cQ,\nabla)$, where $\cQ\in\ECA_{\cO_X}$ and
$\nabla$ is a connection on $\cQ$ give rise in the obvious
way to a stack in Picard groupoids which we denote
$\ECA\nabla_{\cO_X}$.

The ``zero'' object in $\ECA\nabla_{\cO_X}$ is ``the Courant
algebroid'', $\cQ_0 = \Omega^1_X\oplus\cT_X$ with the obvious
connection, the symmetric pairing given by the duality pairing
between $\Omega^1_X$ and $\cT_X$, and the bracket characterized by
the fact the connection is flat.

Note that the pair $(\cQ,\nabla)$ has no non-trivial
automorphisms. The assignment $(\cQ,\nabla)\mapsto c(\nabla)$
gives rise to the morphism of Picard groupoids
\begin{equation}\label{map:curv}
c:\ECA\nabla_{\cO_X} @>>> \Omega^{3,cl}_X \ ,
\end{equation}
where $\Omega^{3,cl}_X$ is viewed as discrete, i.e. the only
morphisms are the identity maps.

\subsubsection{Lemma}
The morphism \eqref{map:curv} is an equivalence.
\begin{pf}
The inverse is given by the following construction. Let
$[\ ,\ ]_0$ denote the Leibniz bracket on $\cQ_0$. For
$H\in\Omega^{3,cl}_X$, $\xi_1,\xi_2\in\cT_X$ let
\[
[\xi_1,\xi_2]_H = [\xi_1,\xi_2]_0 + \iota_{\xi_1\wedge\xi_2}H\ .
\]
This operation extends uniquely to a Leibniz bracket $[\ ,\ ]_H$ on
$\Omega^1_X\oplus\cT_X$ which, together with the symmetric pairing
induced by the duality pairing, is a structure of an exact Courant
algebroid. Let $\cQ_H$ denote this structure.

The obvious connection on  $\cQ_H$ has curvature $H$.
\end{pf}

\subsubsection{Classification of Exact Courant algebroids}
Suppose that $\cQ$ is an exact Courant $\cO_X$-algebroid. The assignment
$\nabla\mapsto c(\nabla)$ gives rise to the morphism
\[
c : \cC(\cQ) @>>> \Omega^{3,cl}_X \ .
\]
The formula $c(\nabla+\alpha) = c(\nabla) + d\alpha$ means that
$c$ is a morphism of sheaves with an action of $\Omega^2_X$, where
$\alpha\in\Omega^2_X$ acts on $H\in\Omega^{3,cl}_X$ by
$H\mapsto H+d\alpha$. Thus, the pair $(\cC(\cQ),c)$ is a torsor under
$(\Omega^2_X @>d>> \Omega^{3cl})$.

\subsubsection{Lemma}
The correspondence $\cQ\mapsto (\cC(\cQ),c)$ establishes an
equivalence of stacks in $\bbC$-vector spaces in categories
\[
\ECA_{\cO_X} @>>> (\Omega^2_X @>d>> \Omega^{3cl})-tors \ .
\]
In particular, the $\bbC$-vector space of isomorphism classes of
exact Courant algebroids is canonically isomorphic to
$H^1(X;\Omega^2_X @>>> \Omega^{3cl})$.

\subsubsection{Locally trivial exact Courant algebroids} An exact
Courant algebroid is said to be {\em locally trivial} if it admits
a flat connection locally on $X$.

Let $\ECA^{loc. triv.}_{\cO_X}$ denote the stack of locally
trivial exact Courant $\cO_X$-algebroids.

\subsubsection{}
For an exact Courant $\cO_X$-algebroid let $\cC^{fl}(\cQ)$ denote
the sheaf of flat connections on $\cQ$.

The sheaf $\cC^{fl}(\cQ)$ is locally nonempty if and only if $\cQ$
is locally trivial, in which case it is a torsor under
$\Omega^{2cl}_X$. The correspondence $\cQ @>>> \cC^{fl}(\cQ)$
establishes an equivalence of stacks $\ECA^{loc. triv.}_{\cO_X}
@>>> \Omega^{2cl}_X-tors$.

\section{Courant extensions of Lie algebroids}
Suppose that $\cA$ a Lie $\cO_X$-algebroid.

\subsection{Courant extensions}

\subsubsection{Definition}
A Courant extension of $\cA$ is a Courant $\cO_X$-algebroid
$\widehat\cA$ together with the isomorphism
$\overline{\widehat\cA}=\cA$ of Lie $\cO_X$-algebroids.

\subsubsection{Morphisms of Courant extensions}
A morphism of Courant extensions of $\cA$ is a morphism
of Courant $\cO_X$-algebroids which is compatible with the
identifications.

Let $\CEXT_{\cO_X}(\cA)$ denote the stack of Courant
extensions of $\cA$.

\subsubsection{}
For a Lie (respectively Courant) $\cO_X$-algebroid $\cA$ let ${\frak g}(\cA)$
denote the kernel of the anchor map. Then, ${\frak g}(\cA)$ is a Lie
(respectively Courant) $\cO_X$-algebroid with the trivial anchor natural in $\cA$.

If $\widehat\cA$ is a Courant extension of $\cA$, then
${\frak g}(\widehat\cA)$ is a Courant extension of ${\frak g}(\cA)$
Hence, there is a morphism ${\frak g}(\ ) : \CEXT_{\cO_X}(\cA) @>>>
\CEXT_{\cO_X}({\frak g}(\cA))$.

\subsection{Courant extensions of transitive Lie
algebroids}\label{subsection:CEXT}
From now on we assume that
$\cA$ is a {\em transitive} Lie $\cO_X$-algebroid and ${\frak
g}(\cA)$ (equivalently $\cA$) is a locally free $\cO_X$-module of
finite rank.

For a sheaf $\cF$ of $\cO_X$-modules let $\dual\cF$ denote the dual
$\shHom_{\cO_X}(\cF,\cO_X)$.

\subsubsection{}
Suppose that $\widehat\cA$ is a Courant extension of $\cA$. Then,
the exact sequence
\begin{equation}\label{ses:Omega-ext}
0 @>>> \Omega^1_X @>>> \widehat\cA @>>> \cA @>>> 0
\end{equation}
is canonically associated to the Courant extension $\widehat\cA$
of $\cA$. Since a morphism of Courant extensions of $\cA$ induces
a morphism of the associated extensions of $\cA$ by $\Omega^1_X$,
it is an isomorphism of the underlying sheaves, and it is clear
that the inverse isomorphism is a morphism of Courant
$\Omega$-extensions of $\cA$.

Therefore, $\CEXT_{\cO_X}(\cA)$ is a stack in groupoids.

\subsubsection{Remark}
$\CEXT_{\cO_X}(\cT_X)$ is none other than $\ECA_{\cO_X}$.

\subsubsection{}\label{sssection:cartesian-square}
Suppose that $\widehat\cA$ is a Courant extension of $\cA$. Let
${\frak g} = {\frak g}(\cA)$ and
$\widehat{{\frak g}} = {\frak g}(\widehat\cA)$ for short.

Since $\langle\widehat{{\frak g}},\Omega^1_X\rangle = 0$, the pairing
on $\widehat\cA$ induces the pairing
\[
\ip : \widehat{{\frak g}}\otimes_{\cO_X}\cA @>>> \cO_X
\]
and  the pairing
\[
\ip : {\frak g}\otimes_{\cO_X}{\frak g} @>>> \cO_X \ .
\]
These yield, respectively, the maps $\widehat{{\frak g}} @>>> \dual\cA$
and ${\frak g} @>>> \dual{{\frak g}}$. Together with the projection
$\widehat{{\frak g}} @>>> {\frak g}$ and  the map $\dual\cA @>>> \dual{{\frak g}}$
adjoint to the inclusion ${\frak g} @>>> \cA$ they fit into the diagram
\begin{equation}\label{diag:pull-back}
\begin{CD}
\widehat{{\frak g}} @>>> \dual\cA \\
@VVV                    @VVV  \\
{\frak g} @>>> \dual{{\frak g}}
\end{CD}
\end{equation}

\subsubsection{Lemma}
The diagram \eqref{diag:pull-back} is Cartesian.

\subsubsection{Corollary}\label{cor:fib-prod}
$\widehat{{\frak g}}$ is canonically isomorphic
to $\dual\cA\times_{\dual{\frak g}}{\frak g}$.

\subsection{Central extensions of Lie
algebras}\label{subsection:c-ext-of-lie}
We maintain the notations
introduced above, i.e. $\cA$ is a transitive Lie $\cO_X$-algebroid
and ${\frak g}$ denotes ${\frak g}(\cA)$ so that there is an exact
sequence
\[
0 @>>> {\frak g} @>i>> \cA @>\pi>> \cT_X @>>> 0 \ .
\]
Hence, $\frak g$ is a Lie algebra in $\cO_X$-modules. We assume
that $\frak g$ is locally free of finite rank.

Suppose in addition that $\frak g$ is equipped with a symmetric $\cO_X$-bilinear
pairing
\[
\ip : {\frak g}\otimes_{\cO_X}{\frak g} @>>> \cO_X
\]
which is invariant under the adjoint action of $\cA$, i.e., for $a\in\cA$ and
$b,c\in{\frak g}$
\[
\pi(a)(\langle b,c\rangle) = \langle [a,b],c\rangle + \langle b,[a,c]\rangle
\]
holds.

\subsubsection{From Lie to Leibniz}\label{sssec:LtoL}
The map $i : {\frak g} @>>> \cA$ and the
pairing on $\frak g$ give rise to the maps
\[
\dual\cA @>\dual i>> \dual{\frak g} @<{\ip}<< {\frak g} \ .
\]
Let $\widehat{\frak g} = \dual\cA\times_{\dual{\frak g}}{\frak g}$ and let
$\pr : \widehat{\frak g} @>>> {\frak g}$ denote the canonical projection.
A section of $\widehat{\frak g}$ is a pair $(\dual a,b)$, where $\dual a\in\dual\cA$
and $b\in\frak g$, which satisfies $\dual i(\dual a)(c) = \langle b,c\rangle$ for
$c\in\frak g$.

The Lie algebra $\frak g$ acts on $\cA$ (by the restriction of the adjoint
action) by $\cO_X$-linear endomorphsims and the map $i : {\frak g} @>>> \cA$
is a map of $\frak g$-modules. Therefore, $\dual\cA$ and $\dual{\frak g}$ are
$\frak g$-modules in a natural way and the map $\dual i$ is a morphism of such.
Hence, $\widehat{\frak g}$ is a $\frak g$-module in a natural way and the map
$\pr$ is a morphism of $\frak g$-modules.

As a consequence, $\widehat{\frak g}$ acquires the canonical structure of a
Leibniz algebra with the Leibniz bracket $[\widehat a,\widehat b]$ of two sections
$\widehat a,\widehat b\in\widehat{\frak g}$ given by the formula
 $[\widehat a,\widehat b] = \pr(\widehat a)(\widehat b)$.

We define a symmetric $\cO_X$-bilinear pairing
\[
\ip : \widehat{\frak g}\otimes_{\cO_X}\widehat{\frak g} @>>> \cO_X
\]
as the composition of $\pr\otimes\pr$ with the pairing on $\frak g$.

The inclusion $\Omega^1_X @>\dual\pi>> \dual\cA$ gives rise to the derivation
$\partial : \cO_X @>>> \widehat{\frak g}$.

\subsubsection{Lemma}
The Leibniz bracket, the symmetric pairing and the derivation defined above
endow $\widehat{\frak g}$ with the structure of a Courant $\Omega^1_X$-extension
of $\frak g$ (in particular, a Courant $\cO_X$-algebroid with the
trivial anchor map).

\subsubsection{Lemma}
The isomorphism of Corollary \ref{cor:fib-prod} is an isomorphism
of Courant extensions of ${\frak g}(\cA)$.

\subsubsection{}\label{sssec:cocycle}
Suppose that $\nabla$ is a connection on $\cA$. $\nabla$ determines
\begin{enumerate}
\item the isomorphism ${\frak g}\bigoplus\cT_X @>\cong>> \cA$
by $a+\xi \mapsto i(a) + \nabla(\xi)$, where $a\in{\frak g}$ and
$\xi\in\cT_X$;
\item the isomorphism $\phi_\nabla : \widehat{{\frak g}} @>>>
\Omega^1_X\bigoplus{\frak g}$ by $(\dual a,b)\mapsto
\dual\nabla(\dual a) + b$, where $\dual a\in\dual\cA$, $b\in\frak g$,
$\dual i(\dual a) = \langle b,\bullet\rangle$,
$\dual\nabla : \dual\cA @>>> \Omega^1_X$ is the transpose of $\nabla$
and $\dual i : \dual\cA @>>> \dual\frak g$ is the transpose of $i$.
\end{enumerate}

Let $[\ ,\ ]_\nabla$ denote the Leibniz bracket on $\Omega^1_X\bigoplus{\frak g}$
induced by $\phi_\nabla$. A simple calculation (which is left to the reader)
shows that it is the extension of the Lie bracket on $\frak g$
by the $\Omega^1_X$-valued (Leibniz) cocycle
$c_\nabla(\ ,\ ) : {\frak g}\otimes{\frak g} @>>> \Omega^1_X$ determined
by $\iota_\xi c(a,b) = \langle [a,\nabla(\xi)],b\rangle$
(where $\xi\in\cT_X$, $a,b\in\frak g$, and the bracket is computed in $\cA$).
In other words, $c_\nabla(\ ,\ )$ is the composition
\[
{\frak g}\otimes{\frak g} @>{\nabla^{\frak g}\otimes\id}>>
\Omega^1_X\otimes{\frak g}\otimes{\frak g} @>{\id\otimes\langle\
,\ \rangle}>> \Omega^1_X \ ,
\]
where $\nabla^{\frak g}$ is the connection on $\frak g$
induced by $\nabla$.

Suppose that $A\in\Omega^1_X\otimes_{\cO_X}{\frak g}$. The automorphism of
$\Omega^1_X\bigoplus{\frak g}$ given by
$(\alpha, a)\mapsto(\alpha+\langle A,a\rangle,a)$ is the isomorphism
of Courant algebroids $(\Omega^1_X\bigoplus{\frak g},[\ ,\ ]_\nabla) @>>>
(\Omega^1_X\bigoplus{\frak g},[\ ,\ ]_{\nabla+A})$ which corresponds to the
identity map on $\widehat{\frak g}$ under the identifications $\phi_\nabla$
and $\phi_{\nabla + A}$.

\subsection{The action of $\ECA_{\cO_X}$}
As before, $\cA$ is a transitive Lie $\cO_X$-algebroid locally free of finite
rank over $\cO_X$, ${\frak g}$ denotes ${\frak g}(\cA)$, $\ip$
is an $\cO_X$-bilinear symmetric $\cA$-invariant pairing on $\frak g$,
$\widehat{{\frak g}}$ is the Courant extension of $\frak g$ constructed in
\ref{sssec:LtoL}.

\subsubsection{}
Let $\CEXT_{\cO_X}(\cA)_\ip$ denote the stack of Courant extensions of $\cA$
which induce the given pairing $\ip$ on $\frak g$. Clearly,
$\CEXT_{\cO_X}(\cA)_\ip$ is a stack in groupoids.

Note that, if $\widehat\cA$ is in $\CEXT_{\cO_X}(\cA)_\ip$, then
${\frak g}(\widehat\cA)$ is canonically isomorphic to $\widehat{{\frak g}}$.

\subsubsection{}
Suppose that $\cQ$ is an exact Courant $\cO_X$ algebroid and
$\widehat\cA$ is a Courant extension of $\cA$. Let $\widehat\cA +
\cQ$ denote the push-out of $\widehat\cA\times_{\cT_X}\cQ$ by the
addition map $\Omega_X^1\times\Omega_X^1 @>+>> \Omega^1_X$. Thus,
a section of $\widehat\cA + \cQ$ is represented by a pair $(a,q)$
with $a\in\widehat\cA$ and $q\in\cQ$ satisfying $\pi(a) =
\pi(q)\in\cT_X$. Two pairs as above are equivalent if their
(componentwise) difference is of the form $(\alpha,-\alpha)$ for
some $\alpha\in\Omega^1_X$.

For $a_i\in\widehat\cA$, $q_i\in\cQ$ with $\pi(a_i)=\pi(q_i)$ let
\begin{equation}\label{formulas:sum}
[(a_1,q_1),(a_2,q_2)] = ([a_1,a_2],[q_1,q_2]),\ \ \
\langle(a_1,q_1),(a_2,q_2)\rangle = \langle a_1,a_2\rangle + \langle q_1,q_2\rangle
\end{equation}
These operations are easily seen to descend to $\widehat\cA + \cQ$. Note that
the compositions
\[
\Omega^1_X @>>> \widehat\cA @>>> \widehat\cA\times_{\cT_X}\cQ
@>>> \widehat\cA + \cQ
\]
and
\[
\Omega^1_X @>>> \cQ @>>> \widehat\cA\times_{\cT_X}\cQ
@>>> \widehat\cA + \cQ
\]
coincide; we denote their common value by
\begin{equation}\label{map:der-sum}
i : \Omega^1_X @>>> \widehat\cA + \cQ \ .
\end{equation}

\subsubsection{Lemma}
The formulas \eqref{formulas:sum} and the map \eqref{map:der-sum}
determine a structure of Courant extension of $\cA$ on
$\widehat\cA + \cQ$. Moreover, the map
${\frak g}(\widehat\cA) @>>> \widehat\cA+\cQ$ defined by $a\mapsto(a,0)$
induces an isomorphim ${\frak g}(\widehat\cA+\cQ)\isomo{\frak g}(\widehat\cA)$
of Courant extensions of ${\frak g}(\cA)$ (by $\Omega^1_X$).

\subsubsection{Lemma}\label{lemma:ECA-action-trans}
Suppose that $\widehat\cA^{(1)}$, $\widehat\cA^{(2)}$ are in
$\CEXT_{\cO_X}(\cA)_\ip$. Then, there exists a unique $\cQ$ in $\ECA_{\cO_X}$,
such that $\widehat\cA^{(2)}= \widehat\cA^{(1)} + \cQ$.
\begin{pf}
Let $\cQ$ denote the quotient of
$\widehat\cA^{(2)}\times_{\cA}\widehat\cA^{(1)}$ by the diagonally
embedded copy of $\widehat{{\frak g}}$. Then, $\cQ$ is an
extension of $\cT$ by $\Omega^1_X$. There is a unique structure of
an exact Courant algebroid on $\cQ$ such that
$\widehat\cA^{(2)}=\widehat\cA^{(1)} + \cQ$.
\end{pf}

\subsection{Courant extensions of flat connections}\label{subsection:CEXT-flat}
Suppose that $\cA$ is a transitive Lie $\cO_X$-algebroid, locally free over $\cO_X$,
${\frak g} = {\frak g}(\cA)$, $\ip = \ip_{\frak g}$ is an $\cA$-invariant
symmetric pairing on $\frak g$.

\subsubsection{}
Suppose that $\nabla$ is a {\em flat} connection on $\cA$ and $\cQ$ is an exact Courant
algebroid. For $a,b\in {\frak g}$, $q,q_1,q_2\in\cQ$ let
\begin{equation}\label{conn-pairing}
\langle a,b\rangle_{\nabla,\cQ} = \langle a,b\rangle_{\frak g},\ \
\langle q_1,q_2\rangle_{\nabla,\cQ} = \langle q_1,q_2\rangle_\cQ,\ \
\langle a,q\rangle_{\nabla,\cQ} = 0\ ,
\end{equation}
and
\begin{equation}\label{conn-bracket}
[a,b]_{\nabla,\cQ} = [a,b]_\nabla,\ \ ,[q_1,q_2]_{\nabla,\cQ} = [q_1,q_2]_\cQ,\ \
[q,a]_{\nabla,\cQ} = \nabla^{\frak g}(\pi(q))(a)
\end{equation}
taking values in ${\frak g}\oplus\cQ$, where $[\ ,\ ]_\nabla$ is the Leibniz bracket
on $\Omega^1_X\oplus{\frak g}$ as in \ref{sssec:cocycle}, $\nabla^{\frak g}$ is
the connection on $\frak g$ induced by $\nabla$ and $\pi$ is the anchor of $\cQ$.

\subsubsection{Lemma}
The formulas \eqref{conn-pairing} and \eqref{conn-bracket} define a structure of a
Courant extension of $\cA$ on ${\frak g}\oplus\cQ$.

\subsubsection{Corollary}
Suppose that $\cA$ admits a flat connection locally on $X$. Then,
$\CEXT_{\cO_X}(\cA)_\ip$ is locally non-empty, hence a torsor
under $\ECA_{\cO_X}$.

\subsubsection{Notation}
We denote ${\frak g}\oplus\cQ$ together with the Courant algebroid structure given by
\eqref{conn-pairing} and \eqref{conn-bracket} by $\widehat\cA_{\nabla,\cQ}$.

\subsubsection{}
Conversely, suppose that $\widehat\cA$ is a Courant extension of $\cA$, and
$\nabla$ is a {\em flat} connection on $\cA$.

Let $\cQ_{\nabla,\widehat\cA}\subset\widehat\cA$ denote the pre-image of $\nabla(\cT_X)$ under the
projection $\widehat\cA @>>> \cA$.

\subsubsection{Lemma}
\begin{enumerate}
\item The (restrictions of the) Leibniz bracket and the symmetric pairing
on $\widehat\cA$ endow $\cQ_{\nabla,\widehat\cA}$ with a structure of an exact
Courant algebroid.

\item $\cQ_{\nabla,\widehat\cA}^\perp\cap\cQ_\nabla = 0$.

\item The projection $\widehat\cA @>>> \cA$ restricts to an isomorphism
$\cQ_{\nabla,\widehat\cA}^\perp @>>> {\frak g}$.

\item $\widehat{\frak g} = \cQ_{\nabla,\widehat\cA}^\perp + \Omega_{\widehat\cA}$.

\item The induced isomorphism $\widehat{\frak g}\isomo {\frak g}\oplus\Omega^1_X$
coincides with the one constructed in \ref{sssec:cocycle}.
\end{enumerate}

\subsubsection{Lemma}\label{lemma:iso-conn}
The isomorphism $\widehat\cA\isomo{\frak g}\oplus\cQ_{\nabla,\widehat\cA}$
induced by $\nabla$ is an isomorphism
$\widehat\cA\isomo\widehat\cA_{\nabla,\cQ_{\nabla,\widehat\cA}}$
(of Courant extensions of $\cA$).

\subsubsection{Change of connection}
Suppose that $\nabla$ is a flat connection on $\cA$ and
$A\in\Omega^1_X\otimes_{\cO_X}{\frak g}$ satisfies the Maurer--Cartan
equation $\nabla A + \displaystyle\frac12[A,A] = 0$, so that the
connection $\nabla + A$ is also flat. Suppose that
$\cQ$ is an exact Courant algebroid. To simplify notations,
let $\widehat\cA = \widehat\cA_{\nabla,\cQ}$.

By Lemma \ref{lemma:iso-conn}, we have the isomorphism
$\widehat\cA\isomo\widehat\cA_{\nabla + A,\cQ_{\nabla+A,\widehat\cA}}$.
Recall that a closed $3$-form $H\in\Omega^{3,cl}_X$ defines an exact Courant
algebroid $\cQ_H$ equipped with a connection, whose curvature is equal to $H$.

Since $A\in\Omega^1_X\otimes_{\cO_X}{\frak g}$ satisfies the Maurer--Cartan
equation, the $3$-form $\langle A,[A,A]\rangle$ is
closed.

\subsubsection{Lemma}
In the notations introduced above, there is an isomorphism
$\cQ_{\nabla + A,\widehat\cA}\isomo \cQ + \cQ_{\langle A,[A,A]\rangle}$.
\begin{pf}
We identify $\cA$ with ${\frak g}\oplus\cT_X$ using the flat connection $\nabla$.
In terms of this identification, the image of $\cT$ under the connection
$\nabla + A$ consists of pairs $(A(\xi),\xi)$, where $\xi\in\cT_X$. Therefore,
$\cQ_{\nabla + A,\widehat\cA}\subset{\frak g}\oplus\cQ=\widehat\cA$ consists
of pairs $(A(\xi),q)$, where $\xi\in\cT_X$ and $\pi(q)=\xi$.

For $i=1,2$, $\xi_i\in\cT_X$, $q_i\in\cQ$ satisfying $\pi(q_i) = \xi_i$ we calculate
the Leibniz bracket in $\widehat\cA$ using \eqref{conn-pairing} and \eqref{conn-bracket}:
\begin{multline*}
[(A(\xi_1),q_1),(A(\xi_2),q_2)]_{\nabla,\cQ} = \\
(A([\xi_1,\xi_2]),
\langle [A(\xi_1),\nabla(\bullet)]_{\frak g},A(\xi_2)\rangle + [q_1,q_2]_\cQ)= \\
(A([\xi_1,\xi_2]), \iota_{\xi_1\wedge\xi_2}\langle [A,A],A\rangle
+ [q_1,q_2]_\cQ)
\end{multline*}
The latter formula shows that the assignment $(A(\pi(q)),q)\mapsto (q,\pi(q))$
viewed as a morphism of extensions of $\cT_X$ by $\Omega^1_X$
\[
\cQ_{\nabla + A,\widehat\cA} @>>> \cQ + (\Omega^1_X\oplus\cT_X)
\]
is, in fact a morphism of exact Courant algebroids
\[
\cQ_{\nabla + A,\widehat\cA} @>>> \cQ + \cQ_{\langle A,[A,A]\rangle} \ .
\]
\end{pf}

\subsubsection{Proposition}\label{prop:loctriv}
For $\cA$ as above, $\nabla$ a flat connection on $\cA$ the
assignment $\cQ\mapsto\widehat\cA_{\nabla,\cQ}$ gives rise to a
morphism $\ECA_{\cO_X} @>>> \CEXT_{\cO_X}(\cA)_\ip$ of
$\ECA_{\cO_X}$-torsors, i.e. a trivialization of
$\CEXT_{\cO_X}(\cA)_\ip$.
\begin{pf}
This amounts to showing that, for $i=1,2$, $\cQ_i\in\ECA_{\cO_X}$,
there is a canonical isomorphism
$\widehat\cA_{\nabla,\cQ_1+\cQ_2}\isomo
\widehat\cA_{\nabla,\cQ_1}+\cQ_2$ which possesses associativity
properties. This follows from \eqref{conn-pairing} and
\eqref{conn-bracket} and the definition of the
$\ECA_{\cO_X}$-action. We leave details to the reader.
\end{pf}

\subsubsection{Corollary}
There is an isomorphism
$\widehat\cA_{\nabla + A,\cQ}\isomo\widehat\cA_{\nabla,\cQ} +
\cQ_{\langle A,[A,A]\rangle}$.

\subsection{The Pontryagin class}\label{ssection:pont}
Suppose that $\cA$ is a transitive Lie algebroid, locally free of
finite rank over $\cO_X$, ${\frak g} = {\frak g}(\cA)$, and
$\langle\ ,\ \rangle$ is an $\cA$-invariant symmetric pairing on
$\frak g$. {\em We assume from now on that $\cA$ admits a flat
connection locally on $X$. }

\subsubsection{}
Suppose that $\cU = \{U_i\}_{i\in I}$ is a covering of $X$ by open subsets
such that $\cA\vert_{U_i}$ admits a flat connection. Let $\nabla_i$ denote
a flat connection on $\cA\vert_{U_i}$. Let $A_{ij} = \nabla_j - \nabla_i$,
$A_{ij}\in\Omega^1_X\otimes_{\cO_X}{\frak g}(U_i\cap U_j)$.

Let $H_{ij} = \langle A_{ij},[A_{ij},A_{ij}]\rangle$,
$H_{ij}\in\Omega^{3,cl}_X(U_i\cap U_j)$; let $H\in\check C^1(\cU;\Omega^{3,cl}_X)$
denote the corresponding cochain.

Let $B_{ijk} =  -\langle A_{ij}\wedge A_{jk}\rangle - \langle A_{jk}\wedge A_{ki}\rangle
+ \langle A_{ki}\wedge A_{ij}\rangle$, $B_{ijk}\in\Omega^2_X(U_i\cap U_j\cap U_k)$;
let $B\in\check C^2(\cU;\Omega^2_X)$ denote the corresponding cochain.

\subsubsection{Lemma}
\begin{enumerate}
\item $\check\partial B = dH = 0$, $\check\partial H = dB$,
i.e. $H+B$ is a $2$-cocycle in $\check C^\bullet(\cU;\Omega^2_X @>>>\Omega^{3,cl}_X)$.

\item The class of $H+B$ in $H^2(X;\Omega^2_X @>>>\Omega^{3,cl}_X)$ does not depend
on the choices of the open cover and of the locally defined flat connections.
\end{enumerate}

\subsubsection{Notation}
We denote the class of $H+B$ in $H^2(X;\Omega^{2}_X@>>>\Omega^{3,cl}_X)$
by $\pont(\cA,\ip)$.

\subsection{The Pontryagin class via obstruction theory}
\subsubsection{}\label{sssection:charcl}
We begin by a brief outline of the cohomological classification of
$\ECA_{\cO_X}$-torsors.

Suppose that $\cS$ is a torsor under $\ECA_{\cO_X}$. Let $\cU = \{U_i\}_{i\in I}$
be a covering of $X$ by open subsets such that $\cS(U_i)$ is non-empty for all $i\in I$.
In this case there are isomorphisms $\ECA_{\cO_X}\vert_{U_i}\isomo\cS\vert_{U_i}$ of
$\ECA_{\cO_X}\vert_{U_i}$-torsors. A choice of such gives rise to the objects
$\cQ_{ij}\in\ECA_{\cO_X}(U_i\cap U_j)$, isomorphims $\cQ_{ij} = - \cQ_{ji}$ and flat
connections $\nabla^0_{ijk}$ on $\cQ_{ijk}=\cQ_{ij}+\cQ_{jk}+\cQ_{ki}$.

We may assume (by passing to a refinement of $\cU$), that the exact Courant algebroids
$\cQ_{ij}$ admit connections. Let $\nabla_{ij}$ denote a connection on $\cQ_{ij}$ whose
curvature we denote by $H_{ij}\in\Omega^{3,cl}_X(U_i\cap U_j)$. Let
$H\in\check C^1(\cU;\Omega^{3,cl}_X)$ denote the corresponding cochain.

The connections $\nabla_{ij}$ give rise to connections on the exact Courant algebroids
$\cQ_{ijk}$. Hence there are $2$-forms $B_{ijk}\in\Omega^2_X(U_i\cap U_j\cap U_k)$ such
that $\nabla_{ijk} = \nabla^0_{ijk} + B_{ijk}$. Let $B\in\check C^2(\cU;\Omega^2_X)$
denote the corresponding cochain.

It is clear that $dH = \check\partial B = 0$ and $\check\partial H = dB$, i.e.
$H+B$ is a $2$-cocycle in $\check C^\bullet(\cU;\Omega^2_X @>>> \Omega^{3,cl}_X)$.
One checks easily that the class of $H+B$ in $H^2(X;\Omega^2_X @>>> \Omega^{3,cl}_X)$
is independent of the choices of the convering $\cU$ and the connections $\nabla_{ij}$.
Moreover, the above construction establishes a bijection (in fact an isomorphism of
$\bbC$-vector spaces) between the isomorphism classes of torsors under $\ECA_{\cO_X}$
and $H^2(X;\Omega^2_X @>>> \Omega^{3,cl}_X)$.

\subsubsection{Theorem}\label{thm:CEXT}
Suppose that $\cA$ is a transitive Lie $\cO_X$-algebroid, locally free of finite rank
over $\cO_X$ which admits a flat connection locally on $X$, and $\ip$
is an $\cA$-invariant symmetric pairing on ${\frak g}(\cA)$
The isomorphism class of the $\ECA_{\cO_X}$-torsor $\CEXT_{\cO_X}(\cA)_\ip$ corresponds
to $\pont(\cA,\ip)$. In particular, Courant extensions of $\cA$ which induce $\ip$ on
${\frak g}(\cA)$ exist globally on $X$ if and only if $\pont(\cA,\ip)=0$.

\begin{pf}
Let $\cU = \{U_i\}_{i\in I}$ be a covering of $X$ by open subsets such that
$\cA\vert_{U_i}$ admits a flat connection. Let $\nabla_i$ be a flat connection
on $\cA\vert_{U_i}$. By Proposition \ref{prop:loctriv}, these give rise
to trivializations of $\CEXT_{\cO_X}(\cA)_\ip\vert_{U_i}$.
The procedure outlined in \ref{sssection:charcl} applied to these yields
the definition of $\pont(\cA,\ip)$ given in \ref{ssection:pont}.
\end{pf}

\subsection{The Pontryagin class Atiyah style}
Throughout this section $\cA$ is a transitive Lie
$\cO_X$-algebroid locally free of finite rank over $\cO_X$, $\frak
g$ denotes ${\frak g}(\cA)$, $\ip$ is a symmetric $\cO_X$-bilinear
$\cA$-invariant pairing on $\frak g$, $\widehat{{\frak g}}$ is the
Courant extension of $\frak g$ as in \ref{sssec:LtoL}.

\subsubsection{The Atiyah class}
The (isomorphism class of the) extension
\begin{equation}\label{ses:algd}
0 @>>> {\frak g} @>>> \cA @>>> \cT_X @>>> 0
\end{equation}
is an element of $\Ext_{\cO_X}(\cT_X, {\frak g})$, whose image under the canonical
isomorphism
$\Ext_{\cO_X}(\cT_X, {\frak g}) @>\cong>> H^1(X;\Omega^1_X\otimes_{\cO_X}{\frak g})$
is called {\em the Atiyah class of $\cA$} and will be denoted $\Atiyah(\cA)$.

Recall that {\em a connection on $\cA$} is a splitting of the extension \eqref{ses:algd}.
Let $\cC(\cA)$ denote the sheaf of locally defined connections on $\cA$. As the difference
of two connections is a map $\cT_X @>>> {\frak g}$, the sheaf $\cC(\cA)$ is a torsor under
$\Omega^1_X\otimes_{\cO_X}{\frak g}$. The Atiyah class $\Atiyah(\cA)$ is the isomorphism
class of the $\Omega^1_X\otimes_{\cO_X}{\frak g}$-torsor $\cC(\cA)$.

The cup product together with the pairing $\ip$ give rise to the map
\[
H^1(X;\Omega^1_X\otimes_{\cO_X}{\frak g})\otimes H^1(X;\Omega^1_X\otimes_{\cO_X}{\frak g})
@>>> H^2(X;\Omega^1_X\otimes_{\cO_X}\Omega^1_X) \ .
\]
We will denote the image of $a\otimes b$
under this map by $\langle a\smile b\rangle$.

\subsubsection{Lemma}
$\langle\Atiyah(\cA)\smile\Atiyah(\cA)\rangle$ is the image of $\pont(\cA,\ip)$
under the map
$H^2(X;\Omega^2_X @>>> \Omega^{3,cl}_X) @>>> H^2(X;\Omega^1_X\otimes_{\cO_X}\Omega^1_X)$.

\subsection{Obstruction theoretic interpretation}
Recall that there is an exact sequence
\[
0 @>>> \Omega^1_X @>>> \widehat{{\frak g}} @>>> {\frak g} @>>> 0 \ .
\]
Spliced with \eqref{ses:algd} it gives rise to the extension
\begin{equation}\label{ext2}
0 @>>> \Omega^1_X @>>> \widehat{{\frak g}} @>>> \cA @>>> \cT_X @>>> 0
\end{equation}
whose isomorphism class is an element of $\Ext^2_{\cO_X}(\cT_X,\Omega^1_X)$.
Let $\beta = \beta(\cA,\ip)$ denote its image under the canonical isomorphism
$\Ext^2_{\cO_X}(\cT_X,\Omega^1_X) @>\cong>> H^2(X;\Omega^1_X\otimes_{\cO_X}\Omega^1_X)$.

The extension \eqref{ext2} gives rise to the stack in groupoids $\cL$
defined as follows. For an open set $U\subseteq X$, $\cL(U)$ is the category of
pairs $(\widehat\cC,\phi)$, where $\widehat\cC$ is a $\Omega^1_X\otimes_{\cO_X}
\widehat{{\frak g}}$-torsor on $U$ and $\phi$ is a morphism of $\Omega^1_X
\otimes_{\cO_X}{\frak g}$-torsors $\widehat\cC\times_{\Omega^1_X\otimes_{\cO_X}
\widehat{{\frak g}}}\Omega^1_X\otimes_{\cO_X}{\frak g} @>>> \cC(\cA)\vert_U$.
Equivalenly, $\cL(U)$ is the category of pairs $(\widehat\cA,\psi)$, where
$\widehat\cA$ is an extension of $\cT_X$ by $\widehat{{\frak g}}$ and $\psi$
is a morphism of the push-out of $\widehat\cA$ by the map $\widehat{{\frak g}}
@>>> {\frak g}$ to $\cA$ (of extensions of $\cT_X$ by $\frak g$).

As is well-known, the stack $\cL$ is a gerbe with lien $\Omega^1_X\otimes_{\cO_X}
\Omega^1$, whose isomorphism class in $H^2(X;\Omega^1_X\otimes_{\cO_X}\Omega^1)$
is, on the one hand, $\beta$, and, on the other hand, the image of $\Atiyah(\cA,\ip)$
under the boundary map $H^1(X;\Omega^1_X\otimes_{\cO_X}{\frak g}) @>>>
H^2(X;\Omega^1_X\otimes_{\cO_X}\Omega^1)$ induced by the extension
\[
0 @>>> \Omega^1_X\otimes_{\cO_X}\Omega^1_X @>>>
\Omega^1_X\otimes_{\cO_X}\widehat{{\frak g}} @>>>
\Omega^1_X\otimes_{\cO_X}{\frak g} @>>> 0 \ .
\]

\subsubsection{Lemma}\label{lemma:c2}
$\beta = \langle\Atiyah(\cA)\smile\Atiyah(\cA)\rangle$.

\subsubsection{Remark}
Note that there is a natural morphism of stacks $\CEXT_{\cO_X}(\cA)_\ip @>>> \cL$
(which maps a Courant extension to the underlying $\cO_X$-module).

\section{Vertex algebroids}
\subsection{Vertex algebroids}
A {\em vertex $\cO_X$-algebroid} is a sheaf of $\bbC$-vector
spaces $\cV$ with a pairing
\begin{eqnarray*}
\cO_X\otimes_\bbC\cV & @>>> & \cV \\
f\otimes v & \mapsto & f*v
\end{eqnarray*}
such that $1* v = v$ (i.e. a ``non-associative unital
$\cO_X$-module'') equipped with
\begin{enumerate}
\item
a structure of a Leibniz $\bbC$-algebra $[\ ,\ ] :
\cV\otimes_\bbC\cV @>>> \cV$

\item
a $\bbC$-linear map of Leibniz algebras $\pi : \cV @>>> \cT_X$
(the {\em anchor})
\item
a symmetric $\bbC$-bilinear pairing $\langle\ ,\ \rangle :
\cV\otimes_\bbC\cV @>>> \cO_X$
\item
a $\bbC$-linear map $\partial : \cO_X @>>> \cV$ such that
$\pi\circ\partial = 0$
\end{enumerate}
which satisfy
\begin{eqnarray}
f*(g*v) - (fg)*v & = & - \pi(v)(f)*\partial(g) -
\pi(v)(g)*\partial(f)\label{assoc} \\
\left[v_1,f*v_2\right] & = & \pi(v_1)(f)*v_2 + f*[v_1,v_2] \label{leib}
\\
\left[v_1,v_2\right] + [v_2,v_1] & = & \partial(\langle v_1,v_2\rangle)
\label{symm-bracket}\\
\pi(f*v) & = & f\pi(v) \label{anchor-lin} \\
\langle f*v_1, v_2\rangle & = & f\langle v_1,v_2\rangle -
\pi(v_1)(\pi(v_2)(f)) \label{pairing}\\
\pi(v)(\langle v_1, v_2\rangle) & = & \langle[v,v_1],v_2\rangle +
\langle v_1,[v,v_2]\rangle \label{pairing-inv} \\
\partial(fg) & = & f*\partial(g) + g*\partial(f) \label{deriv} \\
\left[v,\partial(f)\right] & = & \partial(\pi(v)(f)) \label{bracket-o}\\
\langle v,\partial(f)\rangle & = & \pi(v)(f)\label{pairing-o}
\end{eqnarray}
for $v,v_1,v_2\in\cV$, $f,g\in\cO_X$.

\subsubsection{}
A morphism of vertex $\cO_X$-algebroids is a $\bbC$-linear map of
sheaves which preserves all of the structures.

\subsubsection{Remark}
The notions ``a vertex $\cO_X$ algebroid with the trivial anchor map''
and ``a Courant $\cO_X$ algebroid with the trivial anchor map'' are
equivalent.

\subsection{The associated Lie algebroid}
Suppose that $\cV$ is a vertex $\cO_X$-algebroid. Let
\begin{eqnarray*}
\Omega_\cV & \stackrel{def}{=} & \cO_X*\partial(\cO_X)\subset\cV \ ,\\
\overline\cV & \stackrel{def}{=} & \cV/\Omega_\cV \ .
\end{eqnarray*}

It is easy to see (cf. \cite{B}) that the action of $\cO_X$ on $\cV$
descends to a structure of an $\cO_X$-module on $\overline\cV$ and
the Leibniz bracket on $\cV$ descends to a Lie bracket on $\overline\cV$.
Moreover, there is an evident map $\overline\cV @>>> \cT_X$.

\subsubsection{Lemma}
The $\cO_X$-module $\overline\cV$ with the bracket and the anchor
as above is a Lie $\cO_X$-algebroid.

\subsection{Transitive vertex algebroids}
A vertex $\cO_X$-algebroid is called {\em transitive} if the anchor
map is surjective.

\subsection{Exact vertex algebroids}
A vertex algebroid $\cV$ is called {\em exact} if the map
$\overline\cV @>>> \cT_X$ is an isomorphism. We denote the stack
of exact vertex $\cO_X$-algebroids by $\EVA_{\cO_X}$.

A morphism of exact vertex algebroids induces a morphism of underlying
extensions of $\cT_X$ by $\Omega^1_X$. The latter is an isomorphism and
it is clear that the inverse morphism of extensions is a morphism of vertex
algebroids. Hence, $\EVA_{\cO_X}$ is a stack in groupoids.

\subsubsection{}
It was shown in \cite{B} that $\EVA_{\cO_X}$ is locally non-empty and has
a canonical structure of a torsor under $\ECA_{\cO_X}$. Isomorphism classes
of $\ECA_{\cO_X}$-torsors form a vector space naturally isomorphic to
$H^2(X;\Omega^2_X @>>> \Omega^{3,cl}_X)$. The purpose of this section is
the determination of the isomorphism class of $\EVA_{\cO_X}$. The following
theorem was originally proven in \cite{GMS} by explicit calculations with
representing cocycles. Our proof, ``coordinate-free'' and based on
Theorem \ref{thm:CEXT} and the strategy proposed in \cite{BD}, appears in
\ref{proof:EVA}.

\subsubsection{Theorem (\cite{GMS})}\label{thm:EVA}
The class of $\EVA_{\cO_X}$ in $H^2(X;\Omega^2_X @>>> \Omega^{3,cl}_X)$
is equal to $2\ch_2(\Omega^1_X)$.

\subsubsection{Remark}
Suppose that $P$ is a $GL_n$-torsor on $X$. Let $\cA_P$ denote the Atiyah
algebra of $P$. Then, ${\frak g}(\cA_P) = {\frak gl}_n^P$ and the symmetric
pairing on the latter given by the trace of the product of matrices is
$\cA_P$-invariant. The corresponding Pontriagin class is equal to $2\ch_2(P)$.

\subsection{Vertex extensions of Lie algebroids}
Suppose that $\cA$ is a Lie $\cO_X$-algebroid.

\subsubsection{Definition}
A vertex extension of $\cA$ is a vertex algebroid $\widehat\cA$ together with
an isomorphism $\overline{\widehat\cA}=\cA$ of Lie $\cO_X$-algebroids.

\subsubsection{Morphisms of vertex extensions}
A morphism of vertex extensions of $\cA$ is a morphism of vertex algebroids which
is compatible with the identifications.

Let $\VEXT_{\cO_X}(\cA)$ denote the stack of Courant extensions of
$\cA$.

\subsection{Vertex extensions of transitive Lie algebroids}
From now on we suppose that $\cA$ is a transitive Lie
$\cO_X$-algebroid locally free of finite rank over $\cO_X$. Let
$\widehat{\frak g} = \widehat{\frak g}(\cA)$

\subsubsection{}
Suppose that $\widehat\cA$ is a vertex extension of $\cA$. Then, the derivation
$\partial : \cO_X @>>> \widehat\cA$ induces an isomorphism $\Omega^1_X @>\cong>>
\Omega_{\widehat\cA}$. The resulting exact sequence
\[
0 @>>> \Omega^1_X @>>> \widehat\cA @>>> \cA @>>> 0
\]
is canonically associated to the vertex extension $\widehat\cA$ of $\cA$. Since
a morphism of vertex extensions of $\cA$ induces a morphism of associated extensions
of $\cA$ by $\Omega^1_X$ it is an isomorphism of the underlying sheaves. It is clear
that the inverse isomorphism is a morphism of vertex extensions of $\cA$.

Therefore, $\VEXT_{\cO_X}(\cA)$ is a stack in groupoids.

\subsubsection{Remark}
$\VEXT_{\cO_X}(\cT_X)$ is none other than $\EVA_{\cO_X}$.

\subsubsection{}
Suppose that $\widehat\cA$ is a vertex extension of $\cA$. Let
$\widehat{\frak g} = \widehat{\frak g}(\widehat\cA)$ denote the
kernel of the anchor map (of $\widehat\cA$). Thus, $\widehat{\frak
g}$ is a vertex (equivalently, Courant) extension of $\frak g$.

Analysis similar to that of \ref{subsection:CEXT} shows that
\begin{itemize}
\item the symmetric pairing on $\widehat\cA$ induces a symmetric
$\cO_X$-bilinear pairing on $\frak g$ which is $\cA$-invariant;

\item the vertex extension $\widehat{\frak g}$ is obtained from the
Lie algebroid $\cA$ and the symmetric $\cA$-invariant pairing on
$\frak g$ as in \ref{subsection:c-ext-of-lie}.
\end{itemize}

\subsection{The action of $\ECA_{\cO_X}$}
As before, $\cA$ is a transitive Lie $\cO_X$-algebroid locally
free of finite rank over $\cO_X$, ${\frak g}$ denotes ${\frak
g}(\cA)$, $\ip$ is an $\cO_X$-bilinear symmetric $\cA$-invariant
pairing on $\frak g$, $\widehat{{\frak g}}$ is the Courant
extension of $\frak g$ constructed in \ref{sssec:LtoL}.

\subsubsection{}
Let $\VEXT_{\cO_X}(\cA)_\ip$ denote the stack of Courant
extensions of $\cA$ which induce the given pairing $\ip$ on $\frak
g$. Clearly, $\VEXT_{\cO_X}(\cA)_\ip$ is a stack in groupoids.

Note that, if $\widehat\cA$ is in $\VEXT_{\cO_X}(\cA)_\ip$, then
${\frak g}(\widehat\cA)$ is canonically isomorphic to
$\widehat{{\frak g}}$.

\subsubsection{}
Suppose that $\cQ$ is an exact Courant $\cO_X$ algebroid and
$\widehat\cA$ is a vertex extension of $\cA$. Let $\widehat\cA +
\cQ$ denote the push-out of $\widehat\cA\times_{\cT_X}\cQ$ by the
addition map $\Omega_X^1\times\Omega_X^1 @>+>> \Omega^1_X$. Thus,
a section of $\widehat\cA + \cQ$ is represented by a pair $(a,q)$
with $a\in\widehat\cA$ and $q\in\cQ$ satisfying $\pi(a) =
\pi(q)\in\cT_X$. Two pairs as above are equivalent if their
(componentwise) difference is of the form $(i(\alpha),-i(\alpha))$ for
some $\alpha\in\Omega^1_X$.

For $a\in\widehat\cA$, $q\in\cQ$ with $\pi(a)=\pi(q)$, $f\in\cO_X$
let
\begin{equation}\label{formulas:vert-mult}
f*(a,q) = (f*a, fq),\ \ \ \partial(f) = \partial_{\widehat\cA}(f)
+ \partial_{\cQ}(f)\ .
\end{equation}

For $a_i\in\widehat\cA$, $q_i\in\cQ$ with $\pi(a_i)=\pi(q_i)$ let
\begin{equation}\label{formulas:vert-sum}
[(a_1,q_1),(a_2,q_2)] = ([a_1,a_2],[q_1,q_2]),\ \ \
\langle(a_1,q_1),(a_2,q_2)\rangle = \langle a_1,a_2\rangle +
\langle q_1,q_2\rangle
\end{equation}
These operations are easily seen to descend to $\widehat\cA +
\cQ$.

The two maps $\Omega^1_X @>>> \widehat\cA + \cQ$ given by
$\alpha\mapsto (i(\alpha),0)$ and $\alpha\mapsto (0,i(\alpha))$
coincide; we denote their common value by
\begin{equation}\label{map:der-sum-vert}
i : \Omega^1_X @>>> \widehat\cA + \cQ \ .
\end{equation}

\subsubsection{Lemma}
The formulas \eqref{formulas:vert-mult}, \eqref{formulas:vert-sum}
and the map \eqref{map:der-sum-vert} determine a structure of vertex
extension of $\cA$ on $\widehat\cA + \cQ$. Moreover, the map
the map
${\frak g}(\widehat\cA) @>>> \widehat\cA+\cQ$ defined by $a\mapsto(a,0)$
induces an isomorphim ${\frak g}(\widehat\cA+\cQ)\isomo{\frak g}(\widehat\cA)$
of vertex (equivalently, Courant) extensions of ${\frak g}(\cA)$ (by $\Omega^1_X$).

\subsubsection{Lemma}\label{lemma:ECA-action-trans-vert}
Suppose that $\widehat\cA^{(1)}$, $\widehat\cA^{(2)}$ are in
$\VEXT_{\cO_X}(\cA)_\ip$. Then, there exists a unique $\cQ$ in
$\ECA_{\cO_X}$, such that $\widehat\cA^{(2)}= \widehat\cA^{(1)} +
\cQ$.
\begin{pf}
Let $\cQ$ denote the quotient of
$\widehat\cA^{(2)}\times_{\cA}\widehat\cA^{(1)}$ by the diagonally
embedded copy of $\widehat{{\frak g}}$. Then, $\cQ$ is an
extension of $\cT$ by $\Omega^1_X$. There is a unique structure of
an exact Courant algebroid on $\cQ$ such that
$\widehat\cA^{(2)}=\widehat\cA^{(1)} + \cQ$.
\end{pf}

\subsection{Comparison of $\ECA_{\cO_X}$-torsors}
Suppose that $\widehat\cA$ is a vertex extension of the Lie
algebroid $\cA$. Let $\ip$ denote the induced symmetric pairing
on ${\frak g}(\cA)$.

\subsubsection{}
Suppose that $\cV$ is an exact vertex algebroid. Let $\widehat\cA
-\cV$ denote the pushout of $\widehat\cA\times_{\cT_X}\cV$ by the
difference map $\Omega^1_X\times\Omega^1_X @>{-}>> \Omega^1_X$.
Thus, a section of $\widehat\cA -\cV$ is represented by a pair
$(a,v)$ with $a\in\widehat\cA$, $v\in\cV$ satisfying
$\pi(a)=\pi(v)\in\cT_X$. Two pairs as above are equivalent if
their (componentwise) difference is of the form $(\alpha,\alpha)$
for some $\alpha\in\Omega^1_X$.

For $a\in\widehat\cA$, $v\in\cV$ with $\pi(a)=\pi(v)$, $f\in\cO_X$
let
\begin{equation}\label{formulas:vert-mult-diff}
f*(a,v) = (f*a, f*v),\ \ \ \partial(f) = \partial_{\widehat\cA}(f)
- \partial_{\cV}(f)\ .
\end{equation}

For $a_i\in\widehat\cA$, $v_i\in\cV$ with $\pi(a_i)=\pi(v_i)$ let
\begin{equation}\label{formulas:vert-bracket-diff}
[(a_1,v_1),(a_2,v_2)] = ([a_1,a_2],[v_1,v_2]),\ \ \
\langle(a_1,v_1),(a_2,v_2)\rangle = \langle a_1,a_2\rangle -
\langle v_1,v_2\rangle
\end{equation}
These operations are easily seen to descend to $\widehat\cA - \cV$.

The two maps $\Omega^1_X @>>> \widehat\cA - \cV$ given by
$\alpha\mapsto (i(\alpha),0)$ and $\alpha\mapsto (0,-i(\alpha))$
coincide; we denote their common value by
\begin{equation}\label{map:der-diff}
i : \Omega^1_X @>>> \widehat\cA - \cV \ .
\end{equation}

\subsubsection{Lemma}
The formulas \eqref{formulas:vert-mult-diff}, \eqref{formulas:vert-bracket-diff}
together with \eqref{map:der-diff} determine a structure of a Courant
extension of $\cA$ on $\widehat\cA - \cV$. Moreover, the map
${\frak g}(\widehat\cA) @>>> \widehat\cA-\cV$ defined by $a\mapsto(a,0)$
induces an isomorphim ${\frak g}(\widehat\cA-\cV)\isomo{\frak g}(\widehat\cA)$
of Courant extensions of ${\frak g}(\cA)$ (by $\Omega^1_X$).

\subsubsection{}
The assignment $\cV\mapsto\widehat\cA - \cV$ extends to a functor
\begin{equation}\label{functor:EVAtoCEXT}
\widehat\cA-(\bullet) : \EVA_{\cO_X} @>>> \CEXT_{\cO_X}(\cA)_\ip
\end{equation}
which, clearly, anti-commutes with the respective actions of $\ECA_{\cO_X}$ on
$\EVA_{\cO_X}$ and $\CEXT_{\cO_X}(\cA)_\ip$.

\subsubsection{Corollary}\label{corollary:CEXToppEVA}
The functor \eqref{functor:EVAtoCEXT} is an equivalence of stacks in groupoids.
The isomorphism class of the $\ECA_{\cO_X}$-torsors  $\EVA_{\cO_X}$ and
$\CEXT_{\cO_X}(\cA)_\ip$ are opposite as elements of
$H^2(X;\Omega^2_X @>>> \Omega^{3,cl}_X)$.

\subsubsection{Remark}
In fact, the above construction gives rise to the functor
\[
\VEXT_{\cO_X}(\cA)_\ip @>>> \shHom_{\ECA_{\cO_X}}(\EVA_{\cO_X},\CEXT_{\cO_X}(\cA)_\ip)
\]
which is an equivalence.

\subsection{The canonical vertex $\cO_X$-algebroid}
We will show how the construction of \cite{B} leads to the canonical vertex
extention $\widehat\cA_{\Omega^1_X}^{can}$ of $\cA_{\Omega^1_X}$, the Atiyah algebra
of the cotangent sheaf. As a consequence, we obtain the canonical eqivalence
\[
\widehat\cA_{\Omega^1_X}^{can} - (\bullet) : \EVA_{\cO_X} @>>> \CEXT_{\cO_X}(\cA_{\Omega^1_X})\ ,
\]
which anti-commutes with the action of $\ECA_{\cO_X}$.

\subsubsection{The canonical exact vertex $\Omega^\bullet_X$-algebroid}
All of the considerations regarding the Lie, Courant and vertex algebroids apply
in the differential graded setting. Let $X^\sharp$ denote the dg-madifold with
$\cO_{X^\sharp}$ the de Rham complex of $X$. In \cite{B} we showed that there
exists a unique exact vertex (differential graded) $\cO_{X^\sharp}$-algebroid
which we will denote by $\cU$. Thus, there is a short exact sequence
\[
0 @>>> \Omega^1_{X^\sharp} @>>> \cU @>>> \cT_{X^\sharp} @>>> 0
\]
(of complexes of sheaves on $X$) of which we will be interested in the short exact sequence
\[
0 @>>> \Omega^1_{X^\sharp}{}^{(0)} @>>> \cU^{(0)} @>>> \cT_{X^\sharp}{}^{(0)} @>>> 0
\]
of the degree zero constituents. Note that there is a canonical isomrophism
$\Omega^1_{X^\sharp}{}^{(0)}\isomo \Omega^1_X$.

The natural action of $\cT_{X^\sharp}$ on $\cO_{X^\sharp} = \Omega^\bullet_X$ restricts
to the action of $\cT_{X^\sharp}{}^{(0)}$ on $\cO_X$ and $\Omega^1_X$. The action of
$\cT_{X^\sharp}{}^{(0)}$ on $\cO_X$ gives rise to the map
$\cT_{X^\sharp}{}^{(0)} @>>> \cT_X$ which, together with the natural Lie bracket on
$\cT_{X^\sharp}{}^{(0)}$, endows the latter with a structure of a Lie $\cO_X$-algebroid.

The action of $\cT_{X^\sharp}{}^{(0)}$ on $\Omega^1_X$ gives rise to the map
\begin{equation}\label{map:TtoA}
\cT_{X^\sharp}{}^{(0)} @>>> \cA_{\Omega^1_X} \ ,
\end{equation}
where $\cA_{\Omega^1_X}$ denotes the Atiyah algebra of $\Omega^1_X$.

\subsubsection{Lemma}
The map \eqref{map:TtoA} is an isomorphism of Lie $\cO_X$-algebroids.

\subsubsection{}
It follows that there is a short exact sequence
\[
0 @>>> \Omega^1_X @>>> \cU^{(0)} @>>> \cA_{\Omega^1_X} @>>> 0 \ .
\]

\subsubsection{Lemma}
The vertex $\cO_{X^\sharp}$-algebroid structure on $\cU$ restricts to a
structure of a vertex extension of $\cA_{\Omega^1_X}$ on $\cU^{(0)}$.
The induced symmetric pairing on $\shEnd_{\cO_X}(\Omega^1_X)$ is given by
{\em the negative} of the trace of the product of endomorphisms.

\begin{pf}
The first statement is left to the reader.

According to 5.3 of \cite{B}, the $\cO_X$-vertex algebroid $\cU^{(0)}$ is a
quotient of the ($\cO_X$-vertex algebroid) $\widetilde\cU^{(0)}$, where
\[
\widetilde\cU^{(0)} = \Omega^1_X\bigoplus\left(
\Omega^1_X[1]\otimes_\bbC\cT_X[-1]\bigoplus\cO_X\otimes\cT_X\right)
\]
Moreover, under the quotient map, $\Omega^1_X\bigoplus
\Omega^1_X[1]\otimes_\bbC\cT_X[-1]$ (respectively, $\Omega^1_X[-1]\otimes_\bbC\cT_X[1]$)
surjects onto $\widehat{\frak g}(\cU^{(0)})$ (respectively,
${\frak g}(\cA_{\Omega^1_X})=
\shEnd_{\cO_X}(\Omega^1_X)\isomo\Omega^1_X\otimes_{\cO_X}\cT_X$). The symmetric
pairing on the latter is induced by the one on the former given by the formula
\begin{equation}\label{formula:ip}
\langle \beta_1\otimes\xi_1, \beta_2\otimes\xi_2\rangle =
-(\iota_{\xi_1}\beta_2)(\iota_{\xi_2}\beta_1) \ ,
\end{equation}
where $\beta_i\in\Omega^1_X[-1]$ and $\xi_i\in\cT_X[1]$. 

(The formula for the symmetric pairing on $\widetilde\cU$ in 5.3 of
\cite{B} reads
\[
\langle \beta_1\otimes\xi_1,\beta_2\otimes\xi_2\rangle =
-\beta_1\tau(\xi_2)(\tau(\xi_1)(\beta_2)) -
\beta_2\tau(\xi_1)(\tau(\xi_2)(\beta_1)) -
\tau(\xi_1)(\beta_2)\tau(\xi_2)(\beta_1) \ ,
\]
where $\beta_i\in\cO_{X^\sharp}$, $\xi_i\in\widetilde\cT_X$, $\widetilde\cT_X =
\cT_X[1]\bigoplus\cT_X$ and $\tau$ is the canonical action of $\widetilde\cT_X$
on $\cO_{X^\sharp}$ by derivations with $\xi\in\cT_X[1]$ acting by contraction
$\iota_\xi$. If $\beta_i\in\Omega^1_X$, then the first two summands in the formula
are equal to zero for degree reasons.)

Under the canonical isomorphism
$\shEnd_{\cO_X}(\Omega^1_X)\isomo\Omega^1_X\otimes_{\cO_X}\cT_X$ the symmetric
pairing given by \eqref{formula:ip} corresponds to the one given by {\em the negative}
of the trace of the product of endomorphisms.
\end{pf}

\subsubsection{Corollary}\label{corollary:CEXT-ch}
Let $\ip$ denote the symmetric pairing on ${\frak g}(\cA_{\Omega^1_X}) =
\shEnd_{\cO_X}(\Omega^1_X)$ given by {\em the negative}
of the trace of the product of endomorphisms. Then,
the isomorphism class of $\CEXT_{\cO_X}(\cA_{\Omega^1_X})_\ip$ is equal to
$-2\ch_2(\Omega^1_X)$ in $H^2(X;\Omega^2_X @>>> \Omega^{3,cl}_X)$.

\begin{pf}
The Pontryagin class $\pont(\cA_{\Omega^1_X},\ip)$ is equal to
$-2\ch_2(\Omega^1_X)$. The claim follows from Theorem \ref{thm:CEXT}.
\end{pf}

\subsubsection{Proof of \ref{thm:EVA}}\label{proof:EVA}
Follows from \ref{corollary:CEXToppEVA} and \ref{corollary:CEXT-ch}.


\begin{thebibliography}{ABC}
\bibitem[B]{B} P.~Bressler, Vertex algebroids I, preprint.
\bibitem[BD]{BD} A.~Beilinson, V.~Drinfeld, Chiral algebras, preprint.
\bibitem[GMS]{GMS} V.~Gorbunov, F.~Malikov, V.~Schechtman, Gerbes of chiral
differential operators II, preprint.
\bibitem[S]{S} P.~\v Severa, letters to A.~Weinstein.
\end{thebibliography}
\end{document}